\documentclass[superscriptaddress,amsmath,amssymb,amsfonts,twocolumn,aps]{revtex4-2}
\usepackage{natbib}
\usepackage{amsmath}
\usepackage{amssymb}
\usepackage{amsthm}
\usepackage{nicefrac}
\usepackage{graphicx}
\usepackage{bm}
\usepackage{url}
\usepackage{color}

\newcommand{\bmm}[1]{\bm{\mathrm{#1}}}
\newcommand{\tr}[1]{\mathrm{tr}\left[#1\right]}
\begin{document}
\title{Inferring Structure via Duality for Photonic Inverse Design}
\author{Sean Molesky}
\email{sean.molesky@polymtl.ca}
\affiliation{Department of Engineering Physics, Polytechnique Montr\'{e}al, Montr\'{e}al, Qu\'{e}bec H3T 1J4, CAN}
\author{Pengning Chao}
\affiliation{Department of Mathematics, Massachusetts Institute of Technology, Cambridge, Massachusetts 02139, USA}
\author{Alessio Amaolo}
\affiliation{Department of Chemistry, Princeton University, Princeton, New Jersey 08544, USA}
\author{Alejandro~W.~Rodriguez}
\affiliation{Department of Electrical and Computer Engineering, Princeton University, Princeton, New Jersey 08544, USA}
\begin{abstract}
  \noindent
  Led by a result derived from Sion's minimax theorem concerning constraint violation in quadratically constrained quadratic programs (QCQPs) with at least one constraint bounding the possible solution magnitude, we propose a heuristic scheme for photonic inverse design unifying core ideas from adjoint optimization and convex relaxation bounds. 
  Specifically, through a series of alterations to the underlying constraints and objective, the QCQP associated with a given design problem is gradually transformed so that it becomes strongly dual. 
  Once equivalence between primal and dual programs is achieved, a material geometry is inferred from the solution of the modified QCQP.
  This inferred structure, due to the complementary relationship between the dual and primal programs, encodes overarching features of the optimization landscape that are otherwise difficult to synthesize, and provides a means of initializing secondary optimization methods informed by the global problem context. 
  An exploratory implementation of the framework, presented in a partner manuscript, is found to achieve dramatic improvements for the exemplary photonic design task of enhancing the amount of power extracted from a dipole source near the boundary of a structured material region—roughly an order of magnitude compared to randomly initialized adjoint-based topology optimization for areas surpassing $10~\lambda^{2}$. 
\end{abstract}
\maketitle
\noindent
Extrapolating current trends in computing~\cite{tang2021computing,sevilla2022compute,reuther2022ai,milojicic2021future} and nano-fabrication~\cite{ouyang2023ultrafast,wang2024two,yang2025nanofabrication}, the future of inverse design as a core tool of photonic engineering seems all but certain. 
Either explicitly or implicitly, the understanding of device design as an approximate solution to an optimization problem of the form 
\begin{align}
  &\max/\min~~f_{\circ}\left(\bm{e},\bm{s}\right)
  \label{designProb} \\ 
  &\ni~\bm{e}~\text{respects relevant physics}
  \nonumber \\
  &\land~\bm{s}~\text{respects relevant fabrication tolerances},
  \nonumber
\end{align}
with $\bm{e}$ denoting the electromagnetic field, $\bm{s}$ the variable degrees of freedom defining the device, and $f_{\circ}\left(\bm{e},\bm{s}\right)$ a real-valued objective, seems undeniable. 
As long as algorithms for eq.~\eqref{designProb} continue to improve~\cite{chen2022high,xue2023fullwave,pestourie2023physics}, and real-world use cases continue to complexify~\cite{molesky2018inverse,mao2023multi,shekhar2024roadmapping,kuznetsov2024roadmap}, it is difficult to imagine a scenario in which inverse methods do not become one of the predominant paradigms of photonic device conception within the coming decade.
\\ \\
Accepting this framework, it is natural to question what can be said about its mathematical structure. 
Is it possible to derive expressive limits on achievable values of $f_{\circ}$, for some relatively broad class of objectives, by considering certain aspects of Maxwell's equations, or is the freedom of supposing arbitrary material structuring so great that any sufficiently general limit will invariably have little practical value? 
If realistic limits can be derived, do they imply something useful about high-performance designs?
\\ \\
Perhaps surprisingly, notable progress has been made on this topic over the past five years, despite its loose definition. 
Through the recognition that the constraints of a discretized linear differential equation can be equivalently formulated to restate many versions of eq.~\eqref{designProb} as quadratically constrained quadratic programs (QCQPs), Lagrange duality and semi-definite program (SDP) convex relaxations have proven to be highly effective in treating the first question~\cite{angeris2021convex,chao2022physical}. 
Even with only the constraints of total resistive and reactive power balance, the subtle interplay between chosen material properties, device size, wave physics, and attainable objective values for basic scattering quantities is captured to a remarkable extent~\cite{gustafsson2020upper,angeris2021heuristic}, including the onset of well-known asymptotics such as quasi-statics and geometric optics, and the possibility of creating resonances~\cite{molesky2020fundamental}.
More importantly, through the generalization of imposing complex power balance constraints over any localized volume, convex relaxations have been shown to be remarkably predictive of what is possible for a range of technologically relevant electromagnetic processes given total control over material structure~\cite{zhang2021conservationlawbased,capek2021fundamental,molesky2021comm,shim2024fundamental}. 
Limits coming within a factor of unity of known architectures have been found for a variety of basic scattering and antenna objectives~\cite{gustafsson2019maximum,schab2022upper,sakotic2023perfect,abdelrahman2023thin,li2023transmission,nanda2024exploring}, while limits within factors of ten to one hundred have been found for more challenging applications like Raman-scattering, finite bandwidth radiation enhancement and suppression, and bandwidth-integrated cloaking~\cite{shim2019fundamental,chao2023maximum,strekha2024suppressing,amaolo2024physical,strekha2024limitations}.
\\ \\ 
Still, the general characteristics of eq.~\eqref{designProb} are far from fully resolved. 
While the success of convex relaxation techniques~\cite{candes2010power}, and results like those presented in refs.~\cite{zhong2024topological,guo2024passivity,ma2025time}, intimate that even the simplest photonic design problems subsume a deep mathematical substructure, consensus on the overall utility of convex relaxations and limits for device discovery has yet to form. 
In 2023, ref.~\cite{gertler2025many} provided a first description of how SDP relaxations could be used in photonic structural optimization by applying approximate rank penalization, along with a preliminary application to the design of a one-dimensional multi-layer reflector. 
The results of this study were quite promising, but the simplicity of the considered problem explicitly prevented the possibility of clearly outperforming gradient-based adjoint optimization with random initializations.
The more recent ref.~\cite{dalklint2024performance}, albeit following a protocol differing from what was presented in ref.~\cite{gertler2025many}, has now reported that initializing topology optimization with information obtained from SDP matrix solutions did not meaningfully improve performance in the design of either a small electromagnetic mode converter or resonating elastic plate. 
Moreover, a number of basic questions about duality and SDP relaxations remain open.
Pointedly, there is no broadly accepted explanation as to why observed gaps between best bounds and known device architectures are typically small, often within a factor of about ten, or why these gaps tend to shrink as the allowed device size grows~\cite{chao2023maximum,strekha2024suppressing}. 
Beyond the simplest cases in which the full physics of electromagnetics is restricted down to global power balance, the QCQPs that arise in photonic inverse design problems invariably involve multiple indefinite bilinear forms. 
This (presumptively) broadest class QCQPs does not possess any general guarantees regarding duality gap sizes~\footnote{Please contact us if you are aware of relevant references that we have overlooked concerning this point.}. 
\\ \\
Here, we present a number of conceptual points about the structure of eq.~\eqref{designProb} in relation to quadratic objectives and duality relaxations, and then propose a heuristic ``verlan'' scheme for photonic inverse design based on this understanding. 
Namely, supposing the existence of at least one compact constraint (resistive power balance in a passive photonic device), and a non-empty feasible set (the fields of any physical realization), we apply Sion's minimax theorem to derive four results about Lagrange duality as it applies to ``scattering'' QCQPs, c.f. \emph{\S Definitions and Background}. 
The analysis mainly establishes that as the dual solution field approaches the boundary of any compact constraint it must increasingly respect every other imposed constraints: as the total power consumed by a dual ``polarization'' solution field approaches the power it draws from the exciting field, the dual field becomes increasingly physically realistic. 
This finding may, for example, be used to find the maximum allowed violation of any other scattering constraint given the resistive power balance of the dual solution field and knowledge of the maximum and minimum eigenvalues of the constraint in question. 
It can also be used to approximate how much the inclusion of a constraint that was not previously imposed will tighten a performance limit. 
From this foundation, we then propose three protocols that combine to successively transform any scattering QCQP so that it eventually becomes strongly dual.
These protocols form the basis of our proposed heuristic: once the solution of a closely related strongly dual problem is known, it can be used to infer a globally informed feasible point in the original QCQP; starting from this (hopefully) near-optimal initial point a secondary local optimization can be used to determine a design. 
That is, the scheme proposes a ``verlan'' strategy for inverse design---the reversal of the syllables of ``l'enverse'' (``the reverse'') in the verlan argot---of optimizing through a series of convex relaxations defined from the original QCQP in an attempt to find the ``nearest'' strongly dual formulation and thus determine a near-optimal initial point for a secondary method---reversing the importance typically accorded to initialization and updating in photonic inverse design.
Although these procedures effectively act as analogues of the approximate rank penalization concept employed in ref.~\cite{gertler2025many}, and may be translated to SDP relaxations, their mechanics---proceeding in terms of vectors instead of matrices---are distinct, and their core mechanism presumes structure that is not present in all SDPs~\cite{angeris2021heuristic}. 
It is entirely possible that the two methods could be fused to further improve device discovery. 
\\ \\
A first implementation following the lines of this proposal, achieving roughly order of magnitude improvements compared to adjoint topology optimization with random initializations for areas surpassing $10~\lambda^{2}$ in the exemplary photonic design task of enhancing the amount of power extracted from a dipole source located near the boundary of a structurable material region, is explored in the partner article \emph{Bounds as Blueprints: Towards Optimal and Accelerated Photonics Inverse Design}.
All readers are strongly encouraged to view this article in parallel with the present manuscript. 
\\ \\
The remainder of the text is organized as follows. 
In the next section we review the construction of QCQPs and associated dual programs as they relate to quadratic  objectives in photonic (material structuring) inverse design. 
We then discuss the relation by which we infer structure from the dual solution field, and state why additional control is conceptually desirable---why it is desirable to have a means of pushing the QCQP towards strong duality. 
The next portion of the text derives consequences of Sion's minimax theorem~\cite{sion1958general} for scattering QCQPs. 
The more formal style of this section will hopefully facilitate understanding for researchers more familiar with optimization than photonics. 
With this said, the results we prove are physically intuitive, and the central messages of the article will not be lost if this section is treated as an appendix. 
Drawing from this motivating analysis, we then lay out our verlan scheme for photonic inverse design. 
Finally, the text closes with a brief outlook for future applications. 
\subsection{\S QCQPs in Photonic Device Design}
\noindent
The formulation of eq.~\eqref{designProb} as a QCQP can be done in several equivalent ways~\cite{angeris2019computational,molesky2020hierarchical,kuang2020computational,jelinek2020sub}. 
The approach reviewed here follows the perspective of scattering theory~\footnote{The text does not describe the most general means of forming a QCQP from eq.~\eqref{designProb} via scattering theory. 
For example, any fixed design can be taken as a background, vacuum is only the simplest choice.}, described in substantially greater detail in ref.~\cite{chao2022physical}. 
\\ \\
Take $\bm{e}_{i}$ to be some \emph{initial} electromagnetic field satisfying Maxwell's equations in free space, $\bmm{M}_{\circ}\bm{e}_{i} = \bm{0}$. 
In the presence of a scattering operator $\bmm{X}$ imparting the structure of a photonic device (typically the diagonal operator defined by the dielectric profile, but simply generalizable to magnetic and chiral media) Maxwell's equations become $\left(\bmm{M}_{\circ}-\bmm{X}\right)\left(\bm{e}_{i} + \bm{e}_{s}\right) = \bm{0}$, with $\bm{e}_{i} + \bm{e}_{s} = \bm{e}_{t}$ the total electromagnetic field solution and $\bm{e}_{s}$ the scattered field. 
Letting $\bmm{G}_{\circ}$ be the free space Green function for the design region, the inverse of $\bmm{M}_{\circ}$ with outgoing boundary conditions, the previous relation can be rewritten as $\left(\bmm{M}_{\circ}-\bmm{X}\right)\bm{e}_{s} = \bmm{X}\bm{e}_{i}\Rightarrow\left(\bmm{I}-\bmm{G}_{\circ}\bmm{X}\right)\bm{e}_{t} = \bm{e}_{i}$. 
Restricting this equality to the spatial region occupied by the scatterer, represented by the projection operator $\bmm{I}_{\bmm{X}}$, $\bmm{X}$ is invertible, giving the (integral) equation $\bmm{I}_{\bmm{X}}\left(\bmm{X}^{-1}-\bmm{G}_{\circ}\right)\bmm{X}\bm{e}_{t} = \bmm{I}_{\bmm{X}}\bm{e}_{i}$. 
Trivially, this last expression remains true under composition with any \emph{witness} operator $\bmm{Q}$. 
In the case that $\bmm{Q}$ commutes with $\bmm{I}_{\bmm{X}}$, e.g. the collection of all main diagonal one-hot matrices when considering media with a local polarization response in a spatially local basis for photonic material structuring, it therefore follows that $\bmm{I}_{\bmm{X}}\bmm{Q}\left(\bmm{X}^{-1}-\bmm{G}_{\circ}\right)\bmm{X}\bm{e}_{t} = \bmm{I}_{\bmm{X}}\bmm{Q}\bm{e}_{i}$. 
The quantity $\bm{x} = \bmm{X}\bm{e}_{t}$, up to a scalar unit conversion factor, is the polarization current created within the device in response to the initial field $\bm{e}_{i}$, and is thus restricted to the volume of the scatterer. 
Therefore, $\bmm{X}^{-1}$ of the particular design under consideration (whatever it may be) can be freely replaced with the $\bmm{X}_{\bullet}^{-1}$ inverse scattering potential of the design where scattering media occupies the full design region---the design wherein each controllable parameter is set to its upper limit. 
For the same reason, projecting any such relation along the generated current vector $\bmm{I}_{\bmm{X}}$ is also redundant, $\bm{x}^{\dagger}\bmm{I}_{\bmm{X}} = \bm{x}^{\dagger}$. 
(An equivalent construction for non-local scattering potentials is described in the Appendix \emph{\S QCQP Transformation for Non-Local Scattering}.)
Accordingly, regardless of what the actual scattering potential described by some set of design parameters is, the physical constraints of the scattering system can be stated entirely in terms of the initial field $\bm{e}_{i}$, known operators, and the unknown polarization distribution $\bm{x}$:
\begin{equation}
  \left(\forall\bmm{Q}_{j}\right)~\bm{x}^{\dagger}\bmm{Q}_{j}\bm{e}_{i} - \bm{x}^{\dagger}\bmm{Q}_{j}\left(\bmm{X}^{-1}_{\bullet}-\bmm{G}_{\circ}\right)\bm{x} = 0. 
  \label{sctCnt}
\end{equation}
So long as the symmetric and anti-symmetric parts of $\bm{x}^{\dagger}\bmm{Q}\bmm{X}^{-1}_{\bullet}\bm{x}$ are definite forms, as is the case for any media exhibiting approximately local response electromagnetic response, eq.~\eqref{sctCnt} can be relaxed to an inequality in order to treat continuously variable design parameters:
$$
  \Re\left[\bm{x}^{\dagger}\bmm{Q}_{j}\bm{e}_{i} + \bm{x}^{\dagger}\bmm{Q}_{j}\bmm{G}_{\circ}\bm{x}\right] \left(\geq\lor\leq\right) \Re\left[\bm{x}^{\dagger}\bmm{Q}_{j}\bmm{X}^{-1}_{\bullet}\bm{x}\right],
$$ 
with the correct sign of the $\geq\lor\leq$ determined by the  particularities of the design problem~\footnote{Following a procedure similar to ref.~\cite{amaolo2024can}, these equations can also be extended to treat multi-material design scenarios. 
A similar result, in relation to treating continuous media, is also effectively achieved by limiting imposed $\bmm{Q}_{j}$ operators to volumes that are substantially larger than numerical discretization---allowing $\bmm{X}$ to vary at a scale below the imposed physics is equivalent to allowing the possibility of effective media.}. 
For photonic material structuring problems such a relaxation amounts to allowing ``grey-scale'' media with properties between vacuum and the supposed constituent material. 
In other settings, it could be used, for example, to treat potentials with continuously variable strengths, or describe drive signals restricted to a maximum amplitude.  
If the objective $f_{\circ}$ is a quadratic function of fields that are derivable from $\bm{x}$, the QCQP corresponding to eq.~\eqref{designProb} is then found be asserting all constraints of the form given by eq.~\eqref{sctCnt} for acceptable choices of witnesses---all $\bmm{Q}_{j}$ that commute with any conceivable design~\footnote{In practice, the QCQP used to determined bounds is often furthered relaxed by only considering some special physically motivated subset of choices of $\bmm{Q}$, e.g. $\bmm{Q}$ corresponding to the conservation of power.}. 
From this QCQP, eq.~\eqref{designProb} is convexified by performing the Lagrange duality transformation defined by eq.~\eqref{Lform} and eq.~\eqref{dualSCQP}, or by the SDP relaxation as stated in the Appendix \emph{Equivalence with Semi-Definite Relaxation}. 
\\ \\
While the description given here has focused solely on linear media and quadratic functions of the scattering fields, it should be noted that other versions of eq.~\eqref{designProb} can be treated along conceptually similar lines. 
For example, non-linear processes can be placed in QCQP form by either introducing dummy fields to represent intensities~\cite{mohajan2023fundamental}, or by appealing to the sum-of-squares polynomial hierarchy~\cite{chao2025sum}; certain non-quadratic objectives such as Shannon capacity also have convenient convex relaxations and duality transformations when combined with the scattering constraints~\cite{amaolo2024maximum}.
\subsection{\S Inference} 
\noindent
Shifting perspectives, in a similar spirit to ref.~\cite{gustafsson2024modes}, the scattering relation at the heart of the formulation outlined above can be inverted to connect a generic field with a scattering operator that would allow its interpretation as an induced polarization. 
This interdependence provides a direct structural inference relation: given any field $\bm{x}$ and initial field $\bm{e}_{i}$, there is a linear operator $\bmm{X}_{\bm{x}}$---perhaps with highly unphysical properties---such that if $\bmm{X}_{\bm{x}}$ were realized as a scattering potential, then $\bm{x}$ would be the corresponding polarization field. 
In a sense, the central aim of the entire article is to tame this relationship so that it can be used productively. 
\\ \\
Take $\bm{x}$ to be a discretized complex field over the region available for the design of a photonic device, without any a priori relation to the physical context. 
Letting $l,k$ be indices for the discretization, define 
\begin{align}
  \bmm{X}_{\bm{x}~lk} = \frac{\bm{x}_{k}}{\left(\bmm{G}_{\circ}\bm{x}\right)_{k} + \left(\bm{e}_{i}\right)_{k}}~\delta_{lk}.
  \label{dirInfer}
\end{align}
If $\bmm{X}^{-1}_{\bullet}$ in eq.~\eqref{sctCnt} is replaced by $\bmm{X}^{-1}_{\bm{x}}$, then $\bm{x}$ will satisfy every constraint of the scattering constraint set.
Hence, taking $\bmm{X}^{-1}_{\Delta} = \bmm{X}^{-1}_{\bullet} - \bmm{X}^{-1}_{\bm{x}}$, with $\!^{\mathsf{s}}$ denoting symmetrization $\left[\bmm{A}\right]^{\mathsf{s}} = \bmm{A} + \bmm{A}^{\dagger}$, the constraint $\Re\left[\bm{x}^{\dagger}\bmm{Q}_{j}\bm{e}_{i}\right] - \bm{x}^{\dagger}\left[\bmm{Q}_{j}\left(\bmm{X}^{-1}_{\bullet}-\bmm{G}_{\circ}\right)\right]^{\mathsf{s}}\bm{x} = 0$ is recast as
\begin{equation}
  \bm{x}^{\dagger}\left[\bmm{Q}_{j}\bmm{X}^{-1}_{\Delta}\right]^{\mathsf{s}}\bm{x} = 0.
  \label{matCns}
\end{equation}
The content of eqs.~\eqref{dirInfer} and \eqref{matCns} is twofold. 
First, given a field $\bm{x}$, satisfaction of a relation of the form given by eq.~\eqref{sctCnt} can be viewed either as a statement of local complex power balance~\cite{kuang2020computational}, or as a statement that the implied scattering potential obeys the physical material constraint reported by the witness $\bmm{Q}_{j}$. 
Second, by using eq.~\eqref{dirInfer} and projecting to the nearest feasible scattering potential under the matrix norm, every field---and in particular every ``super-optimal'' field computed during the calculation of duality bounds---directly implies a material structure. 
As shown in ref.~\cite{chao2025blueprints}, this simple recipe for relating dual solution fields to structures, which is by no means unique, works astonishingly well.
In fact, it appears that just two specially chosen witnesses (representing global resistive and reactive power balance) may be enough for eq.~\eqref{matCns} to effectively pilot gradient-based adjoint optimization towards greatly improved geometries over large enough design regions. 
However, the ``jump'' from the implicit scattering potential to the inferred design can in principle be very large, potentially leading to an effective disconnect between the limit solution field and the actualized polarization field. 
The analysis of the following two sections is carried out, in brief, towards the aim of developing a reliable procedure for controlling the distance that must be spanned. 
\section{Sion's Theorem and Duality}
\noindent
With photonics in mind, we present implications of Sion's minimax theorem for a subclass of QCQPs, here represented by $\mathtt{P}$, that appears to occur frequently in physical inverse design problems~\cite{angeris2021convex,gertler2025many}. 
Specifically, after making preliminary definitions and showing that Sion's theorem applies (lemma 0) eq.~\eqref{sion}, we prove four results related to $\mathtt{P}$ and its Lagrangian $\mathcal{L}\left(\bm{\phi},\bm{x}\right)$, eq.~\eqref{Lform}.
\\ \\
\emph{Lemma 0.} 
The value of the Sion program $\mathtt{S}\left(\mathtt{P}\right)$, eq.~\eqref{SionFunc}---given by inverting the order of maximization and minimization in the dual program $\mathtt{D}\left(\mathtt{P}\right)$, eq.~\eqref{dualSCQP}---is equal to the value of the dual program $\mathtt{D}\left(\mathtt{P}\right)$. 
\\ \\
\emph{Lemma 1.} 
The Sion set $S_{\mathtt{P}}$, defined below eq.~\eqref{SionFunc}, contains the convex hull $C_{\mathtt{P}}$ of feasible points for $\mathtt{P}$. 
By lemma 4, the two sets agree on the boundary of any constraint with a positive definite bilinear part, and exact equivalence implies strong duality. 
\\ \\
\emph{Lemma 2.} 
Let $\bm{\phi}_{\circledast}$ be a solution of $\mathtt{D}\left(\mathtt{P}\right)$. 
The difference between an $\bm{x}_{\circledast}$ solution of the Sion program $\mathtt{S}\left(\mathtt{P}\right)$, eq.~\eqref{SionFunc}, and the internal optimal $\tilde{\bm{x}}_{\circledast}$ determined when computing $\mathtt{D}\left(\mathtt{P}\right)$ at $\bm{\phi}_{\circledast}$ is confined to the kernel of the bilinear part of $\mathcal{L}_{\bm{\phi}_{\circledast}}\left(\bm{x}\right)$ ($\mathcal{L}\left(\bm{x},\bm{\phi}_{\circledast}\right)$ with $\bm{\phi}_{\circledast}$ held fixed). 
That is, beyond the fact that $\mathtt{D}\left(\mathtt{P}\right)$ and $\mathtt{S}\left(\mathtt{P}\right)$ have equivalent values, lemma 0, they must also determine largely equivalent $\bm{x}$ and $\bm{\phi}$ vectors. 
\\ \\
\emph{Lemma 3.} 
In many relevant cases, the multiplier cone $\Psi_{\mathtt{P}}$ that enters the definitions of $\mathtt{D}\left(\mathtt{P}\right)$ and $\mathtt{S}\left(\mathtt{P}\right)$ can be bounded without altering the solutions of either program.
If this change is made, the Sion function, eq.~\eqref{SionFunc}, becomes continuous. 
\\ \\
\emph{Lemma 4.} 
Approaching the boundary of a constraint with a positive definite bilinear part from within the Sion set $S_{\mathtt{P}}$ implies shrinking violation in all constraints. 
\\ \\
While lemmas 1--3 furnish helpful supporting details, the framework presented in \emph{\S Verlan Scheme} is primarily advised by lemma 4---any QCQP can be shifted towards strong duality by pushing the maximum of the Sion function towards the boundary of a compact constraint. 
Readers willing to accept this result as stated may freely jump to the next section without losing the main points we hope to convey. 
\subsection{\S Definitions and Background}
\noindent
A (complex) \emph{quadratic} function is a rule of the form
$$
  f:\mathbb{C}^{n}\rightarrow\mathbb{R}\ni\bm{x}\mapsto 2\Re\left[\bm{s}^{\dagger}\bm{x}\right] - \bm{x}^{\dagger}\bmm{A}\bm{x}+c,
$$
with $\bmm{A}:\mathbb{C}^{n}\rightarrow\mathbb{C}^{n}$ Hermitian, $\bm{s}^{\dagger}\in\mathbb{C}^{n\dagger}$ and $c\in\mathbb{R}$. 
A quadratic function is \emph{positive definite} or \emph{compact} (resp. \emph{positive semi-definite}) if $\bmm{A}\succ 0$ ($\bmm{A}\succeq 0$.) 
Motivated by eq.~\eqref{sctCnt}, letting an $\mathsf{s}$ superscript denote the symmetrization $\bmm{M}^{\mathsf{s}} = \bmm{M} +\bmm{M}^{\dagger}$, we define a set of quadratic functions $\left\{f_{j}\right\}_{j\in J}$ to be a set of \emph{scattering constraints} if $\left(\forall j\in J\right)$ either (a) $f_{j}\left(x\right)$ is linear, $\bmm{A}_{j} = \bmm{0}$, or 
$$
  \text{(b)}~f_{j}:\bm{x}\mapsto 2\Re\left[\bm{s}^{\dagger}\bmm{Q}_{j}\bm{x}\right] - \bm{x}^{\dagger}\left[\bmm{Q}_{j}\bmm{U}\right]^{\mathsf{s}}\bm{x} + c_{j}
$$
with $\bmm{Q}_{j}$ selected from some vector space of admissible \emph{witness} operators~\footnote{Connecting with \emph{\S Appendix: Non-Local Scattering Constraints}, the witness terminology is based on the notion of ``witnessing'' a certain design parameter, in the scattering form, of the overarching device optimization problem.}, $\bmm{Q}_{j}\in Q$, $\bm{s}\in\mathbb{C}^{n}$ and $\bmm{U}:\mathbb{C}^{n}\rightarrow\mathbb{C}^{n}$ fixed over the indexing collection. 
We define a \emph{scattering constrained quadratic program} (\emph{SCQP}) as an optimization problem $\mathtt{P}$ of the form 
\begin{align}
  &\max_{\bm{x}\in F_{\bm{\kappa}}}f_{o}\left(\bm{x}\right)
  \label{SCQP} \\
  &\ni \left(\forall j\in J\right) f_{j}\left(\bm{x}\right) \geq 0
  \nonumber,
\end{align}
wherein $f_{o}$ is a quadratic objective function, $J$ is some finite indexing for a set of scattering constraints, feasible points exist, and there is a set of scalars $\bm{\kappa} = \left\{\kappa_{j}\right\}_{j\in J}$ such that the non-negative domain of  
$$
  f_{\bm{\kappa}}\left(\bm{x}\right) = \sum_{j\in J}\kappa_{j} f_{j}\left(\bm{x}\right),
$$
$F_{\bm{\kappa}}$, $\left(\forall j\in J\right)\kappa_{j}\geq 0$, is compact~\footnote{It is plausible that the majority of the ideas we have investigated here could be extended to a wider class of objective functions. 
Namely, we suspect that analogous behaviour holds whenever the dual transformation is well defined.}: A SCQP is a feasible QCQP such that (a) the imposed constraints form a scattering constraint set, and (b) some subset of the scattering constraint set can be combined to form a composite compact constraint. 
The assumptions of feasibility and compactness are often obvious in physical device design: any feasible selection of design parameters will lead to a polarization field satisfying every scattering constraints that can be imposed through the underlying differential equations; in well-posed problems there are some guard rails preventing the possibility of unbounded measurables given a bounded input. 
\\ \\
Although seemingly restrictive, SCQPs are sufficiently expressive to state the subset sum, an NP-complete problem widely considered to be computationally difficult to solve exactly~\cite{karp2009reducibility,aaronson2005guest}. 
\\ \\
Subset sum: Let $S\! =\! \left\{s_{1},\ldots, s_{n}\right\}$ be a set of integers and $t\in\mathbb{Z}$. 
Does some subset of $S$ sum to $t$?
\\ \\
Let $\bm{z}\in\mathbb{Z}^{n}\subset\mathbb{C}^{n}$ be given by taking the element of $S$ as the components of a vector, and $\bm{1}\in\mathbb{C}^{n}$ be the vector with a value of $1$ in every entry. 
The subset sum program $\mathtt{P}$ for $\left(S,t\right)$ is then 
\begin{align}
  &\max_{\bm{x}\in F_{\bm{\kappa}}} \Re\left[\bm{z}^{\dagger}\bm{x}\right]
  \label{subSum} \\
  &\ni\left(\forall k\leq n\right)~2\Re\left[\bm{1}^{\dagger}\bmm{\delta}_{kk}\bm{x}\right]-\bm{x}^{\dagger}\left[\bmm{\delta}_{kk}\right]^{\mathsf{s}}\bm{x} = 0
  \nonumber \\
  &\land\left(\forall k\leq n\right)~2\Re\left[\bm{1}^{\dagger}i\bmm{\delta}_{kk}\bm{x}\right]-\bm{x}^{\dagger}\left[i\bmm{\delta}_{kk}\right]^{\mathsf{s}}\bm{x} = 0
  \nonumber \\
  &\land~\Re\left[\bm{1}^{\dagger}\bm{x}\right]\geq 1~\land~\Re\left[\bm{z}^{\dagger}\bm{x}\right]\leq t.
  \nonumber
\end{align}
The first two lists of constraints, which form a scattering set, restrict each entry of $\bm{x}$ to be either zero or one. 
The first linear constraint forces at least one entry of $\bm{x}$ to be non-zero, and the second linear constraint forces the sum to be no larger than $t$ (in order to match our max convention). 
The first list of constraints can be summed to form a compact constraint. 
As long as some subset of $S$ sums to a value less than $t$, which is a condition for $\left(S,t\right)$ to be non-trivial, $\mathtt{P}$ is feasible. 
\\ \\
Any function formed as a sum of the defining scattering quadratics of a SCQP $\mathtt{P}$,
$$
  f_{\bm{\phi}}\left(\bm{x}\right) = 2\Re\left[\bm{s}^{\dagger}\bmm{Q}_{\bm{\phi}}\bm{x}\right] - \bm{x}^{\dagger}\left[\bmm{Q}_{\bm{\phi}}\bmm{U}\right]^{\mathsf{s}}\bm{x} + c_{\bm{\phi}} = \sum_{j\in J}\phi_{j}f_{j}\left(\bm{x}\right), 
$$ 
is a \emph{composite} constraint. 
$F_{\bm{\phi}} = \left\{\bm{x}\in\mathbb{C}^{n}~|~f_{\bm{\phi}}\left(\bm{x}\right)\geq 0\right\}$ is used to denote the feasible set of $f_{\bm{\phi}}$---so that that the $F_{\bm{\kappa}}$ set appearing in eq.~\eqref{SCQP} is compact and convex~\cite{bredon2013topology}. 
\\ \\ 
Setting $\bm{s}_{\bm{\phi}} = \sum_{J}\phi_{j}\bmm{Q}_{j}\bm{s}$, $\bmm{A}_{\bm{\phi}} = \sum_{J}\phi_{j}\left[\bmm{Q}_{j}\bmm{U}\right]^{\mathsf{s}}$, $\bm{s}_{\bm{\psi}} = \bm{s}_{o} + \bm{s}_{\bm{\phi}}$, $\bmm{A}_{\bm{\psi}} = \bmm{A}_{o} + \bmm{A}_{\bm{\phi}}$, and $c_{\bm{\phi}} = \sum_{J}\phi_{j}c_{j}$, the \emph{Lagrangian} of a SCQP $\mathtt{P}$ , like that of a QCQP, is
\begin{align}
  \mathcal{L}\left(\bm{\phi}, \bm{x}\right) =
  2\Re\left[\bm{s}_{\bm{\psi}}^{\dagger}\bm{x}\right]\!-\! 
  \bm{x}^{\dagger}\bmm{A}_{\bm{\psi}}\bm{x} + c_{\bm{\phi}}.
  \label{Lform}
\end{align}
\begin{figure}[t!]
  \vspace{-12pt}
  \centering
  \includegraphics[width=1.0\columnwidth]{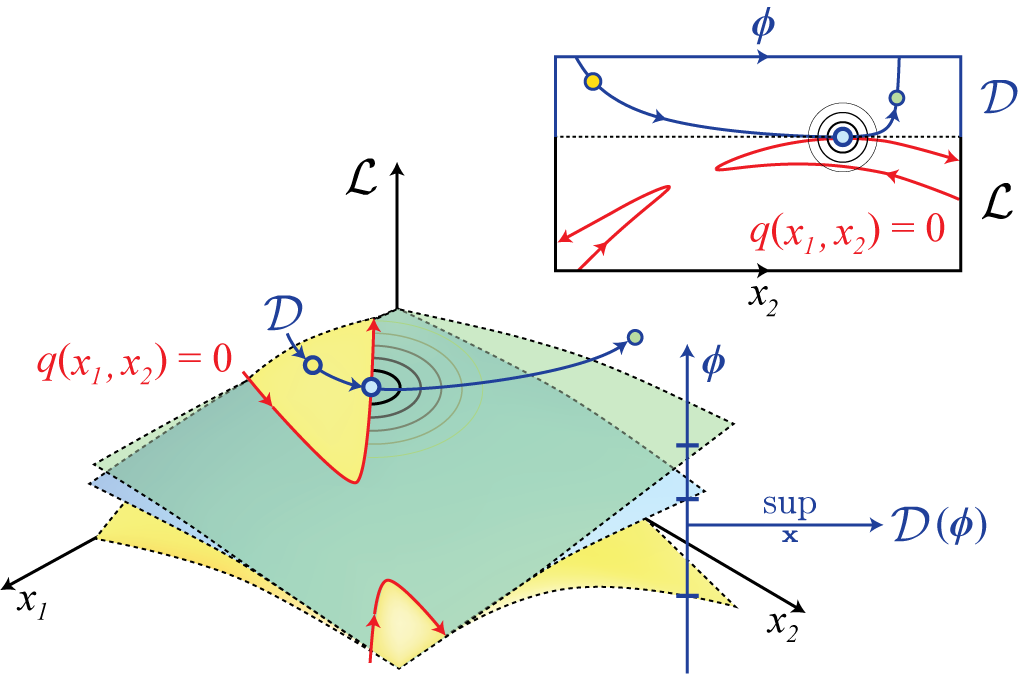}
  \caption{\bmm{Pictorial Lagrangian duality.} 
  The left part of the figure depicts three sections, $\phi$ multiplier values, of a schematic optimization Lagrangian $\mathcal{L}$ for a two dimensional $\left(x_{1},x_{2}\right)$ QCQP subject to a single feasible equality constraint, the red path intersecting the three sections. 
  The value of the (convex) dual $\mathcal{D}$, determined within each section by the maximum of $\mathcal{L}$ over the associated manifold (represented by the coloured balls), is necessarily at least as large as the maximum of $\mathcal{L}$ over the constraint set, the rippling circle. 
  In the illustration, as indicated by the upper inset, the minimum of $\mathcal{D}$, the maximum of the light blue section, is exactly equal to the maximum of $\mathcal{L}$ along the constraint path.
  For a general QCQP no such agreement between the maximum of $\mathcal{L}$ within the set of points satisfying all constraints and the minimum of $\mathcal{D}$ need occur, with the two values generally differing by a duality gap. 
  An optimization program exhibits strong duality when the two values agree.
  }
  \label{dualFig}
\vspace{-12pt}
\end{figure}
The \emph{dual} (program) of $\mathtt{P}$, see fig.~\ref{dualFig}, is given by
\begin{equation}
  \mathtt{D}\left(\mathtt{P}\right)
  = \inf_{\bm{\phi}\in\Psi_{\mathtt{P}}}\max_{\bm{x}\in F_{\bm{\kappa}}}\mathcal{L}\left(\bm{\phi},\bm{x}\right) 
  = \inf_{\bm{\phi}\in\Psi_{\mathtt{P}}}\mathcal{D}\left(\bm{\phi}\right), 
  \label{dualSCQP}
\end{equation} 
with 
$
  \mathcal{D}\left(\bm{\phi}\right) = \max_{\bm{x}\in F_{\bm{\kappa}}}\mathcal{L}\left(\bm{\phi},\bm{x}\right)
$
the dual function, $\Psi_{\mathtt{P}} = \{\bm{\phi}\in\mathbb{H}^{J}~|~\bmm{A}_{\bm{\phi}} \succeq \bmm{A}_{o}\}$, and $\mathbb{H}$ the positive half space $\mathbb{H} = \left\{r\in\mathbb{R}~|~r\geq 0\right\}$---$\Psi_{\mathtt{P}}$ is the closure of $\{\bm{\phi}\in\mathbb{H}^{J}~|~\bmm{A}_{\bm{\phi}}\succ \bmm{A}_{o}\}$. 
For later use we also define $\Phi_{\mathtt{P}} = \{\bm{\phi}\in\mathbb{H}^{J}~|~\bmm{A}_{\bm{\phi}}\succeq\bmm{0}\}$. 
Since the supremum of a collection of convex lower semi-continuous functions is both convex and lower semi-continuous, the dual function of any SCQP is convex and lower semi-continuous. 
A SCQP is \emph{strongly dual} if $\mathtt{D}\left(\mathtt{P}\right) = \mathtt{P}$---the value of the dual program coincides with the value of the \emph{primal} program.
$\circledast$ subscripts are used to denote optimal points. 
For example, $\bmm{\phi}_{\circledast}$ would be used to denote a minimizer of $\mathcal{D}\left(\bm{\phi}\right)$. 
The $\#$ symbol is used to denote an active variable in a partially applied function. 
\\ \\
Note that the definition of $\mathcal{D}\left(\bm{\phi}\right)$ given above is unconventional: typically $\max_{\bm{x}\in F_{\bm{\kappa}}}$ is replaced by $\sup_{\bm{x}\in\mathbb{C}^{n}}$. 
So long as a compact constraint exists, this alteration only produces superficial differences.  
The convexity of $\Psi_{\mathtt{P}}$ implies that it possible to come arbitrarily close to the minimum of $\mathcal{D}\left(\bm{\phi}\right)$ along a path $p:\left[0,1\right]\rightarrow\Psi_{\mathtt{P}}$ such that $\bmm{A}_{p\left(t\right)}\succ \bmm{A}_{o}$. 
The Hessian of $\mathcal{D}\left(\bm{\phi}\right)$ is well defined along $p\left(t\right)$, and if 
$
  \bm{x}_{\circledast} = \bmm{A}^{-1}_{\bm{\psi}}\bm{s}_{\bm{\psi}}\not\in F_{\bm{\kappa}},
$
there is a $\delta > 0$ such that increasing $\bm{\psi}$ along $\bm{\kappa}$ will result in $\bmm{A}_{\bm{\psi}+\delta\bm{\kappa}}\succ\bmm{A}_{\bm{\psi}}$ and $\mathcal{D}\left(\bm{\psi}+\delta\bm{\kappa}\right) < \mathcal{D}\left(\bm{\psi}\right)$~\footnote{A similar argument  shows that $\ker\bmm{A}_{\bm{\phi}_{\circledast}}\neq\bmm{0}$, then $\bm{s}_{\circ}\not\in\ker\bmm{A}_{\bm{\phi}_{\circledast}}$ unless $\inf_{\bm{\phi}\in\Psi_{\mathtt{P}}}\mathcal{D}\left(\bm{\phi}\right) = \mathcal{D}\left(\bm{0}\right)$~\cite{molesky2020hierarchical}.}. 
That is, if the standard dual is minimized with sufficiently fine precision, the restriction $\bm{x}\in F_{\bm{\kappa}}$ is immaterial. 
Imposing $\max_{x\in F_{\bm{\kappa}}}$ is technically helpful as it guarantees the continuity of the dual, which is used in lemma 4, c.f. Appendix. 
\\ \\
The remainder of the article centres on the application of the following theorem to SCQPs.
\\ \\
\emph{Sion's minimax theorem~\cite{sion1958general,komiya1988elementary}.}
Let $F_{\bm{\kappa}}$ be a compact convex subset of a topological vector space, $\Psi_{\mathtt{P}}$ a convex subset of a topological vector space, and $\mathcal{L}$ a real valued function on $\Psi_{\mathtt{P}}\times F_{\bm{\kappa}}$. 
Suppose (a) that $\mathcal{L}\left(\#,\bm{x}\right)$ is lower semicontinuous and quasi-convex on $\Psi_{\mathtt{P}}$ for each $\bm{x}\in F_{\bm{\kappa}}$, and (b) that $\mathcal{L}\left(\bm{\phi},\#\right)$ is upper semicontinuous and quasi-concave on $F_{\bm{\kappa}}$ for each $\bm{\phi}\in\Psi_{\mathtt{P}}$. 
Then,
\begin{equation}
  \max_{\bm{x}\in F_{\bm{\kappa}}}\inf_{\bm{\phi}\in\Psi_{\mathtt{P}}}\mathcal{L}\left(\bm{\phi},\bm{x}\right) = 
  \inf_{\bm{\phi}\in\Psi_{\mathtt{P}}}\max_{\bm{x}\in F_{\bm{\kappa}}}~\!\mathcal{L}\left(\bm{\phi},\bm{x}\right).
  \label{sion} 
\end{equation}
\noindent 
Mirroring the definition of $\mathcal{D}\left(\bm{\phi}\right)$, we take 
\begin{equation}
  \mathcal{S}\left(\bm{x}\right) = \inf_{\bm{\phi}\in\Psi_{\mathtt{P}}}\mathcal{L}\left(\bm{\phi},\bm{x}\right)
  \label{SionFunc}
\end{equation}
to be the \emph{Sion} function of $\mathtt{P}$, and $\mathtt{S}\left(\mathtt{P}\right) = \max_{\bm{x}\in F_{\bm{\kappa}}}\mathcal{S}\left(\bm{x}\right)$ to be the Sion program of $\mathtt{P}$. 
Because the infimum of a collection of concave upper semi-continuous functions retains these properties, the Sion function of any SCQP is quasi-concave and upper semi-continuous. 
\\ \\
In discussing the behaviour of Sion functions we have found it useful to consider the (non-empty) feasible set of a SCQP $F_{\mathtt{P}}$, the convex hull of $F_{\mathtt{P}}$, $C_{\mathtt{P}}$, and the Sion set, $S_{\mathtt{P}} = \left\{\bm{x}\in F_{\bm{\kappa}}~|~\mathcal{S}\left(\bm{x}\right) \neq -\infty\right\}$.
For example, simplifying to the case that $f_{o}$ is simply linear, if $\bm{x}\in F_{\mathtt{P}}$ then $\mathcal{S}\left(\bm{x}\right) = f_{o}\left(\bm{x}\right)$---every constraint is satisfied and so $\bm{\phi} = \bm{0}$ is the unique minimizer of $\mathcal{L}\left(\bm{\phi},\bm{x}\right)$. 
Hence, it follows that $\mathcal{S}\left(\bm{x}\right) = f_{o}\left(\bm{x}\right)$ if $\bm{x}\in C_{\mathtt{P}}$.
\subsection{\S Duality in SCQPs}
\noindent
\emph{Lemma 0.} Eq.~\eqref{sion} applies to eq.~\eqref{SCQP}: 
$$
  \inf_{\bm{\phi}\in\Psi_{\mathtt{P}}}\mathcal{D}\left(\bm{\phi}\right) 
  = \max_{\bm{x}\in F_{\bm{\kappa}}}\mathcal{S}\left(\bm{x}\right).
$$
\begin{proof}
  As $\bm{\phi}\in\Psi_{\mathtt{P}}\Rightarrow\bmm{A}_{\bm{\psi}}\succeq\bmm{0}$, $\mathcal{L}\left(\bm{\phi},\#\right) = \mathcal{L}_{\bm{\phi}}\left(\bm{x}\right)$ is a concave continuous function on $F_{\bm{\kappa}}$, ref.~\cite{folland1999real}, that is bounded from above by $\bm{s}^{\dagger}_{\bm{\psi}}\bmm{A}_{\bm{\psi}}^{-1}\bm{s}_{\bm{\psi}} + c_{\bm{\phi}}$.
  Suppose $\bm{x},\bm{y}\in F_{\bm{\kappa}}$. 
  For any $t\in\left(0,1\right)$ 
  \begin{align}
    &t\mathcal{L}_{\bm{\phi}}\left(\bm{x}\right) + \left(1-t\right)\mathcal{L}_{\bm{\phi}}\left(\bm{y}\right) 
    = 2\Re\left[\bm{s}_{\bm{\psi}}^{\dagger}\left(t\bm{x}+\left(1-t\right)\bm{y}\right)\right] - \nonumber \\ 
    & t~\bm{x}^{\dagger}\bmm{A}_{\bm{\psi}}\bm{x} - \left(1-t\right)
    \bm{y}^{\dagger}\bmm{A}_{\bm{\psi}}\bm{y} + c_{\bm{\phi}} =
    \nonumber \\
    &\mathcal{L}_{\bm{\phi}}\left[t\bm{x}+\left(1-t\right)\bm{y}\right]
    -t\left(1-t\right)\left(\bm{x}-\bm{y}\right)^{\dagger}
    \bmm{A}_{\bm{\psi}}\left(\bm{x}-\bm{y}\right) \Rightarrow
    \nonumber \\
    &\mathcal{L}_{\bm{\phi}}\left[t\bm{x}+\left(1-t\right)\bm{y}\right]
    \geq t\mathcal{L}_{\bm{\phi}}\left(\bm{x}\right) + \left(1-t\right)\mathcal{L}_{\bm{\phi}}\left(\bm{y}\right)~\Rightarrow
    \nonumber \\
    &t\bm{x}+\left(1-t\right)\bm{y}\in F_{\bm{\kappa}}.
    \nonumber
  \end{align}
  $\Psi_{\mathtt{P}}$ convex as $\bmm{A}_{\bm{\theta}}\succeq\bmm{A}_{o}$ and $\bmm{A}_{\bm{\phi}}\succeq\bmm{A}_{o}$ imply
  $$
    \bm{x}^{\dagger}\left[t\bmm{A}_{\bm{\phi}} +\left(1-t\right)\bmm{A}_{\bm{\theta}}\right]\bm{x}\geq\bm{x}^{\dagger}\bmm{A}_{o}\bm{x},
  $$ 
  $\Rightarrow t\bm{\phi} + \left(1-t\right)\bm{\theta}\in\Psi_{\mathtt{P}}$.
  $\mathcal{L}\left(\#,\bm{x}\right) = \mathcal{L}_{\bm{x}}\left(\bm{\phi}\right)$ is convex and continuous on $\Psi_{\mathtt{P}}$ since it is affine. 
\end{proof}
\noindent
\emph{Lemma 1}. If $f_{o}$ is linear, then $\bm{y}\in C_{\mathtt{P}}\Rightarrow\mathcal{S}\left(\bm{y}\right) = f_{o}\left(\bm{y}\right)$.
\\ \\
\emph{Corollary.} $C_{\mathtt{P}}\subseteq S_{\mathtt{P}}$ for every SCQP. 
\begin{proof}
  Let $\left\{\bm{x}_{j}\right\}_{j\in J}\subset F_{\mathtt{P}}\ni$
  $
    \bm{y} = \sum_{j\in J}t_{j}\bm{x}_{j}\in C_{\mathtt{P}},
  $ 
  with $\left(\forall j\in J\right)~0<t_{j}<1$ and $\sum_{j\in J} t_{j} = 1$. 
  Partition $J$ into $J_{a}$ and $J_{b}$, so that $\bm{y}_{a} = \sum_{j\in J_{a}}t_{k}\bm{x}_{j}$, $\bm{y}_{b} = \sum_{j\in J_{b}}t_{j}\bm{x}_{j}$, $n_{a} = \sum_{j\in J_{a}}t_{j}$ and $n_{b} = \sum_{j\in J_{b}}t_{j}$. 
  Since $\mathcal{S}$ is concave 
  $
    \mathcal{S}\left(\bm{y}\right) \geq n_{a}~\mathcal{S}\left(\bm{y}_{a}/n_{a}\right) + n_{b}~\mathcal{S}\left(\bm{y}_{b}/n_{b}\right).
  $
  Taking $J_{a}$ equal to a single element, this result specializes to $\mathcal{S}\left(\bm{y}\right) \geq t_{j}~\mathcal{S}\left(\bm{x}_{j}\right) + n_{b}~\mathcal{S}\left(\bm{y}_{b}/n_{b}\right) = t_{j}f_{o}\left(\bm{x}_{j}\right) + n_{b}~\mathcal{S}\left(\bm{y}_{b}/n_{b}\right)$. 
  Therefore, $\mathcal{S}\left(\bm{y}\right) \geq \sum_{j\in J} t_{j}~\mathcal{S}\left(\bm{x}_{j}\right)$ $\geq\sum_{j\in J} t_{j}f_{o}\left(\bm{x}_{j}\right) = f_{o}\left(\bm{y}\right)$, since $f_{o}$ is assumed to be linear. 
  Because $\mathcal{S}\left(\bm{y}\right)\leq f_{o}\left(\bm{y}\right)$, $\mathcal{S}\left(\bm{y}\right) = f_{o}\left(\bm{y}\right)$.
  \\ \\
  The corollary follows as $\mathcal{S}\left(\bm{y}\right)\geq\sum_{j\in J} t_{j}~\mathcal{S}\left(\bm{x}_{j}\right)\geq\sum_{j\in J} t_{j}f_{o}\left(\bm{x}_{k}\right)$ holds generally. 
\end{proof}
\noindent
\emph{Remark.} The inclusion stated in the corollary is typically strict, and differences between $S_{\mathtt{P}}$ and $C_{\mathtt{P}}$ are a central cause of duality gaps: The maximum of a linear function on $C_{\mathtt{P}}$ is equivalent to its maximum on $F_{\mathtt{P}}$, but its maximum on $S_{\mathtt{P}}$ may be greater. 
\\ \\
\emph{Lemma 2.} Let $\bm{\phi}_{\circledast}$ minimize $\inf_{\bm{\phi}\in\Psi_{\mathtt{P}}}\max_{\bm{x}\in F_{\bm{\kappa}}}\mathcal{L}\left(\bm{\phi},\bm{x}\right)$, and take $\bm{x}_{\circledast} = \bmm{A}^{-1}_{\bm{\psi}_{\circledast}}\bm{s}_{\bm{\psi}_{\circledast}}$---with $\bmm{A}^{-1}_{\bm{\psi}_{\circledast}}$ the pseudo inverse of $\bmm{A}_{\bm{\psi}_{\circledast}}$ in the event that $\bmm{A}_{\bmm{\psi}_{\circledast}}\not\succ 0$. 
If $\tilde{\bm{x}}_{\circledast}$ is a solution of $\max_{\bm{x}\in F_{\bm{\kappa}}}\inf_{\bm{\phi}\in\Psi_{\mathtt{P}}}\mathcal{L}\left(\bm{\phi},\bm{x}\right)$, then there exists $\bmm{k}\in\ker{\bmm{A}_{\bmm{\phi}_{\circledast}}}\ni\tilde{\bm{x}}_{\circledast} = \bm{x}_{\circledast} + \bmm{k}$. 
\\ \\
\emph{Corollary.} If $\bm{x}_{\circledast a}$ and $\bm{x}_{\circledast b}$ are associated with the dual minimizers $\bm{\phi}_{\circledast a}$ and $\bm{\phi}_{\circledast b}$ as above, then $\bm{x}_{\circledast a} = \bm{x}_{\circledast b}$ on the complement of $\ker\bmm{A}_{\bm{\psi}_{\circledast a}}\cup\ker\bmm{A}_{\bm{\psi}_{\circledast b}}$. 
If $\left\{\bm{\phi}_{\circledast j}\right\}_{j\in J}$ is a collection of dual minimizers, then two solutions of $\max_{\bm{x}\in F_{\bm{\kappa}}}\inf_{\bm{\phi}\in\Psi_{\mathtt{P}}}\mathcal{L}\left(\bm{\phi},\bm{x}\right)$ can differ only an element of
$$ 
  \bigcap_{j\in J}\ker\bmm{A}_{\bm{\psi}_{\circledast j}}.
$$
\begin{proof}
  Let $\bm{\phi}_{\circledast}$ be the minimizer of $\inf_{\bm{\phi}\in\Psi_{\mathtt{P}}}\mathcal{D}\left(\bm{\phi}\right)\Rightarrow\bmm{A}_{\bmm{\psi}_{\circledast}}\succeq\bmm{0}$. 
  $\bm{x}_{\circledast} = \bmm{A}^{-1}_{\bm{\psi}_{\circledast}}\bm{s}_{\bm{\psi}_{\circledast}}$ is then a solution of $\max_{\bm{x}\in F_{\bm{\kappa}}}\mathcal{L}\left(\bm{\phi_{\circledast}},\bm{x}\right)$. 
  By lemma 0, there is then an $\tilde{\bm{x}}_{\circledast}\in F_{\bm{\kappa}}$ such that 
  $
    \inf_{\bm{\phi}\in\Psi_{\mathtt{P}}}\mathcal{L}\left(\bm{\phi},\tilde{\bm{x}}_{\circledast}\right) = 
    \mathcal{L}(\bm{\phi}_{\circledast}, \bm{x}_{\circledast}). 
  $
  Because $\bm{x}_{\circledast}$ is selected by $\max_{\bm{x}\in\mathtt{C}}\mathcal{L}\left(\bm{\phi}_{\circledast},\bm{x}\right)$, 
  $$
   \mathcal{L}\left(\bm{\phi}_{\circledast},\bm{x}_{\circledast}\right) = \inf_{\bm{\phi}\in\Psi_{\mathtt{P}}}\mathcal{L}\left(\bm{\phi},\tilde{\bm{x}}_{\circledast}\right)\leq 
   \mathcal{L}(\bm{\phi}_{\circledast},\tilde{\bm{x}}_{\circledast}) \leq 
   \mathcal{L}\left(\bm{\phi}_{\circledast},\bm{x}_{\circledast}\right).
  $$
  Therefore, $\mathcal{L}(\bm{\phi}_{\circledast},\tilde{\bm{x}}_{\circledast}) = \mathcal{L}(\bm{\phi}_{\circledast},\bm{x}_{\circledast})\Rightarrow\tilde{\bm{x}}_{\circledast} = \bm{x}_{\circledast} + \bmm{k}$ with $\bmm{k}\in\ker\bmm{A}_{\bm{\phi}_{\circledast}}$. 
  \\ \\
  The corollary follows as analogous argument applies to under the exchange of order---maximize first. 
\end{proof}
\noindent 
\emph{Lemma 3.} Let $\Phi_{\mathtt{P}}^{b} = \left\{\bm{\phi}\in\mathbb{H}^{J}~|~\bmm{A}_{\bm{\phi}}\succeq\bmm{0}~\land~\lVert\bm{\phi}\rVert \leq b\right\}$ and $\Psi_{\mathtt{P}}^{b} = \left\{\bm{\phi}\in\mathbb{H}^{J}~|~\bmm{A}_{\bm{\phi}}\succeq\bmm{A}_{o}~\land~\lVert\bm{\phi}\rVert \leq b\right\}$. 
If there exists a $\delta > 0$ such that $\forall\bm{\phi}\in\partial\Phi_{\mathtt{P}}^{1}$
$$
  \max_{\bm{x}\in F_{\bm{\kappa}}}f_{\bm{\phi}}\left(\bm{x}\right) > \delta,
$$ 
there is a finite number $b > 0$ such that  
\begin{align}
  \max_{\bm{x}\in F_{\bm{\kappa}}}\min_{\bm{\phi}\in\Psi_{\mathtt{P}}^{b}}\mathcal{L}\left(\bm{\phi},\bm{x}\right)
  = \max_{\bm{x}\in F_{\bm{\kappa}}}\inf_{\bm{\phi}\in\Psi_{\mathtt{P}}}\mathcal{L}\left(\bm{\phi},\bm{x}\right).
  \nonumber
\end{align}
That is, $\max_{x\in C}\mathcal{S}\left(\bm{x}\right) = \max_{x\in C}\mathcal{S}_{b}\left(\bm{x}\right)$ with 
$$
  \mathcal{S}_{b}\left(\bm{x}\right) = \min_{\bm{\phi}\in\Psi_{\mathtt{P}}^{b}}\mathcal{L}\left(\bm{\phi},\bm{x}\right)
$$ 
continuous and concave. 
\begin{proof}
  Appealing to compactness, let $u$ and $l$ be upper and lower bounds of $f_{o}$ over $F_{\bm{\kappa}}$.
  As $\bm{0}\in\Phi_{\mathtt{P}}^{b}$, and each of the multiplier sets is convex, if $\bm{\phi}\in\Psi_{\mathtt{P}}$ then $\bm{\phi}/\lVert\bm{\phi}\rVert\in\Phi_{\mathtt{P}}^{1}$. 
  Set $b = \left(u - l\right)/\delta$. 
  If $\lVert\bm{\phi}\rVert \geq b~\land~\bm{\phi}\in\Psi_{\mathtt{P}}$, then $\max_{x}\mathcal{D}\left(\bm{\phi}\right)\geq u$. 
  Therefore, $\inf_{\bm{\phi}\in\Psi_{\mathtt{P}}}\mathcal{D}\left(\bm{\phi}\right) = \min_{\bm{\phi}\in\Psi^{b}_{\mathtt{P}}}\mathcal{D}\left(\bm{\phi}\right)$. 
  \\ \\
  Because $\Psi_{\mathtt{P}}^{b} = \Psi_{\mathtt{P}}\cap\Phi_{\mathtt{P}}^{b}$ is convex, $\Psi^{b}_{\mathtt{P}}$ can be substituted for $\Psi_{\mathtt{P}}$ in lemma 0 without alteration. 
  Since $\Psi^{b}_{\mathtt{P}}$ is closed and bounded, it is compact, and so there $\exists\bm{\phi}_{\circledast}\in\Psi^{b}_{\mathtt{P}}$ such that $\mathcal{D}\left(\bm{\phi}_{\circledast}\right) = \inf_{\bm{\phi}\in\Psi_{\mathtt{P}}}\mathcal{D}\left(\bm{\phi}\right)$. 
  $\mathcal{S}_{b}\left(\bm{x}\right)$ is continuous by lemma A \emph{Appendix}.  
\end{proof}
\noindent
\emph{Remark.} The conditions of the lemma are satisfied, for example, if $\bm{\phi}\in\partial\Phi^{1}_{\mathtt{P}}\Rightarrow c_{\bm{\phi}} > 0\lor\lVert\bm{s}_{\bm{\phi}}\rVert > 0$. 
This follows from the fact that $\partial\Phi^{1}_{\mathtt{P}}$ is compact, $\bm{\phi}\in\partial\Phi^{1}_{\mathtt{P}}$ implies that the maximum eigenvalue of $\bmm{A}_{\bm{\phi}}$ is bounded, and that $c_{\bm{\phi}}$ and $\lVert\bm{s}_{\bm{\phi}}\rVert$ are continuous functions of $\bm{\phi}$.
\\ \\
\emph{Lemma 4.} If $\bm{x}_{\circledast}$ maximizes $\max_{\bm{x}\in F_{\bm{\kappa}}}\inf_{\bm{\phi}\in\Psi_{\mathtt{P}}}\mathcal{L}\left(\bm{\phi},\bm{x}\right)$, and $f_{\bm{\gamma}}$ is a positive definite equality constraint such that $f_{\bm{\gamma}}\left(\bm{x}_{\circledast}\right) = 0$, then $f_{j}\left(\bm{x}_{\circledast}\right)\geq 0$ for all $j\in J$.
\begin{proof}
  Without loss of generality, suppose that there is a feasible point with objective value greater than zero, and $f_{j}\left(\bm{x}_{\circledast}\right) = -k < 0$. 
  There is then a $\delta >0$ such that $f_{j}\left(\bm{x}_{\circledast}\right) + \delta f_{\bm{\gamma}}\left(\bm{x}_{\circledast}\right) = -k$ and $\bmm{A}_{j} + \delta \bmm{A}_{\bm{\gamma}}\succeq \bmm{A}_{o}$. 
  Let $\bm{\zeta} = \phi_{j} + \delta\bm{\gamma}\in\Psi_{\mathtt{P}}$, and set $m = \Re\left[\bm{s}_{o}^{\dagger}\bm{x}_{\circledast}\right]$. 
  Then, for any $n > \max\left\{m, 1/k\right\}$, $\mathcal{L}\left(n\bm{\zeta}/k,\bm{x}_{\circledast}\right) < 0$, contradicting the assumption that $\mathtt{P}$ possesses feasible solutions with objective values greater than zero. 
\end{proof}
\noindent
\emph{Remark.} The lemma restricts possible constraint violations of an optimal $\bm{x}_{\circledast}$ based on the tightness of any positive definite constraint. 
Simplifying to the case that $\bmm{A}_{o} = \bmm{0}$, take $f_{j}\left(\bm{x}_{\circledast}\right) = -r$ and $f_{\gamma}\left(\bm{x}_{\circledast}\right) = s$, with $\bmm{A}_{\bm{\gamma}}\succ 0$ and $r,s >0$. 
Let $n$ to be the most negative eigenvalue of $\bmm{A}_{j}$ and $m$ to be the smallest eigenvalue of $\bmm{A}_{\gamma}$. 
Because $f_{j} + \left|n/m\right|f_{\gamma}$ is then a positive semi-definite constraint, the definition of $\bm{x}_{\circledast}$ requires that $\left|n/m\right|s\geq r$. 
The same argument also implies that $S_{\mathtt{P}}$ coincides with $F_{\mathtt{P}}$ on $\partial F_{\bm{\kappa}}$. 
\\ \\
Lemma 4 guides the remainder of the manuscript: there is a relation between strong duality, allowable constraint violation, and the nearness of a $\bm{x}_{\circledast}$ maximizer of $\max_{\bm{x}\in F_{\bm{\kappa}}}\inf_{\bm{\phi}\in\Psi_{\mathtt{P}}}\mathcal{L}\left(\bm{\phi},\bm{x}\right)$ to the boundary of a positive definite constraint. 
In particular, if a maximizer exist on the boundary of $F_{\bm{\kappa}}$, strong duality holds. 
\section{Proposal and Outlook}
\begin{figure*}[t!]
 \vspace{-12pt}
 \centering
 \includegraphics[width=2.0\columnwidth]{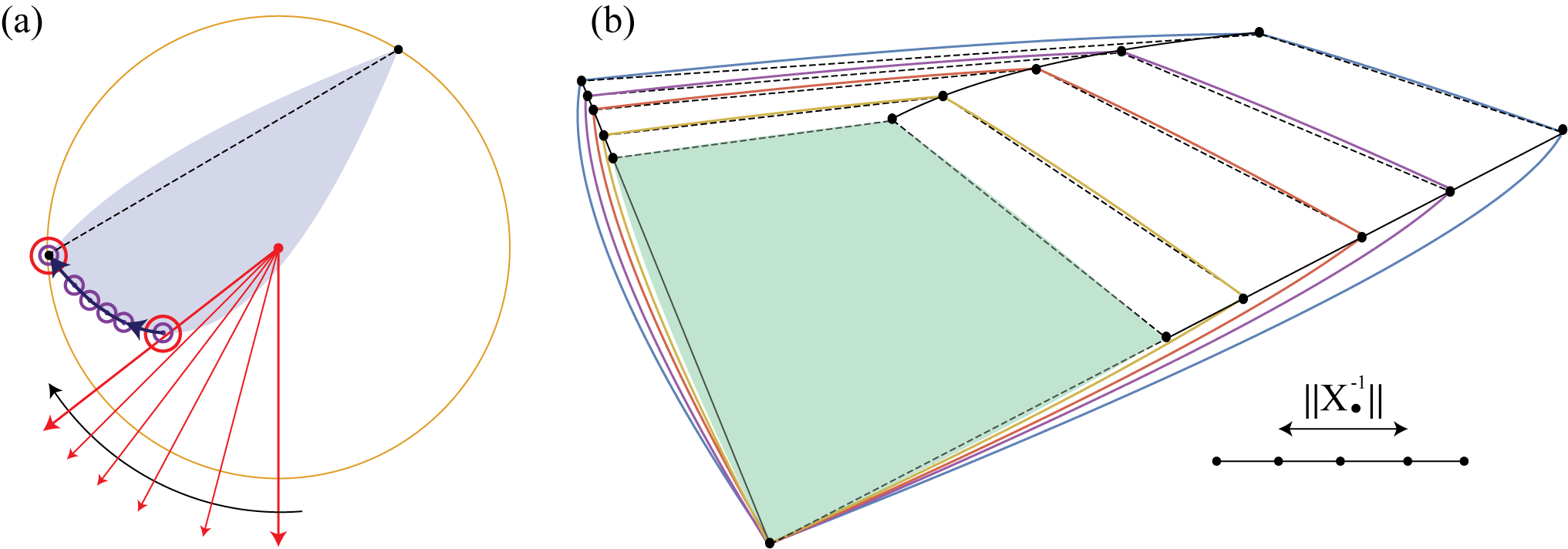}
 \caption{\bmm{Schematic protocols.} 
  The figure illustrates the scrape and contract protocols of our proposed scheme. 
  (a) A series of scraping steps for the program $\max \left(0,-1\right)^{\dagger}\bm{x},~\bm{x}\in\mathbb{R}^{2}$ subject to the indefinite constraint $4 + 4x_{a} -3x_{b} - 4x^{2}_{b} = 0$ and the compact constraint $1 - x^{2}_{a} - x^{2}_{b} = 0$ (orange disk). 
  The tear-drop shaped purple region is the Sion set, bounding the convex hull depicted as a dashed line.  
  Taking the dual transform of $\mathtt{P}$, the minimum is found to be $\bm{x} = \left(0,-1/3\right)$---the point emphasized by the bullseye. 
  Consequently there is a duality gap of $1/3$ compared to $\mathtt{P}$: $\mathtt{P}$ has two feasible points, the intersection points of the Sion set with the orange circle, with the best objective value being achieved at $\bm{x} = \left(0,0\right)$. 
  Smaller purple circles show the location of the $\bm{x}_{\circledast}$ solution of $\mathtt{D}\left(\mathtt{P}\left(\bm{r}\right)\right)$ for five successive scrapes: $\bm{r}^{\left(n\right)} = \bm{r}^{\left(o\right)} / 2 + \bm{x}_{\circledast}^{\left(o\right)}/ 2$ with $\bm{r}$ initialized as $\bm{s}_{o}$. 
  As the objective vector drifts, $\bm{x}_{\circledast}$ is observed to move towards the boundary of the orange disk.
  On the fifth modification of $\bm{s}_{o}$ strong duality is achieved. 
  (b) Shape of the Sion set $S_{\mathtt{P}}$ (boundaries demarcated by solid coloured outlines and feasible points) compared to the convex hull $C_{\mathtt{P}}$ (boundaries demarcated by dashed outlines and feasible points) under contraction for an example scattering constraint set in $\mathbb{R}^{2}$.
  As $\lVert\bmm{X}_{\bullet}^{-1}\rVert$ becomes increasingly positive definite, both $S_{\mathtt{P}}$ and $C_{\mathtt{P}}$ contract towards the trivial feasible point at the origin, the common lowest point of all regions. 
  Simultaneously, $S_{\mathtt{P}}$ increasingly coincides with $C_{\mathtt{P}}$, implying strong duality for any objective with maximum occurring outside of $S_{\mathtt{P}}$.  
  }
  \label{sionFig}
\end{figure*}
\noindent
The results given above are certainly only a small part of what can be understood about Lagrange duality via Sion's theorem. 
Nevertheless, they are already sufficient to devise heuristics for SCQPs that appear to work quite well for photonic inverse design~\cite{chao2025blueprints}.  
The algorithm elements suggested in the next section are based on three observations concerning \emph{\S Duality in SCQPs}. 
\\ \\
(Oa) 
Let $\mathtt{P}$ be a generic SCQP with $\bm{s}_{o}$ the linear part of objective. 
The Sion set $S_{\mathtt{P}}$ of $\mathtt{P}$ is effectively determined by the imposed constraints, with $\bmm{A}_{o}$ acting only as an offset and $\bm{s}_{o}$ playing no role.
Hence, viewing $\bm{s}_{o}$ as variable parameter, there are always instances of $\bm{s}_{o}$ such that $\mathtt{P}\left(\bm{s}_{o}\right)$ is strongly dual.  
\\ \\
(Ob) 
A shrinking of the separation between the dual solution $\bm{x}_{\circledast}$ and the boundary of a compact constraint $\partial F_{\bm{\kappa}}$ implies a shrinking of the duality gap, lemma 4. 
A reasonable criterion for modifying $\bm{s}_{o}$ so that $\mathtt{P}\left(\bm{s}_{o}\right)$ becomes strongly dual is, correspondingly, to decrease the distance between $\bm{x}_{\circledast}$ and $\partial F_{\bm{\kappa}}$. 
C.f. fig.2 (a). 
\\ \\
(Oc) 
In any SCQP resulting from a scattering theory the $\bmm{U}$ operator of a scattering constraint set naturally decomposes into variable $\bmm{X}^{-1}_{\bullet}$ and fixed $-\bmm{G}_{\circ}$ terms---the variable material structure controlled by the designer and the background physics. 
If $\left(\bmm{X}^{-1}_{\bullet}-\bmm{G}_{\circ}\right)\bm{x} = \bm{s}$ on the domain where $\bm{x}$ is non-zero, then, so long as $\bm{x}$ satisfies all explicitly linear constraints and $\left(\forall j\in J\right)~c_{j}\leq 0$, $\bm{x}$ will be feasible.  
This relation of feasibility is continuous as a function of $\bmm{X}^{-1}_{\bullet}$---in varying $\bmm{X}^{-1}_{\bullet}$ the difference vectors between new and old feasible points can be restricted to arbitrarily small volumes by controlling the matrix norm of updates. 
If $\bmm{X}^{-1}_{\bullet}$ is made increasingly positive definite, as witnessed by some composite constraint, the solution set of the SCQP will contract, causing $S_{\mathtt{P}}$ to approach $C_{\mathtt{P}}$. 
For instance, as $\bmm{X}^{-1}_{\bullet}\rightarrow\alpha\bmm{I}$, with $\alpha\rightarrow\infty$, the feasible set of a SCQP collapses towards the zero vector---in the absence of a scattering medium, the only feasible solution is the trivial solution.
C.f. fig.2 (b). 
\\ \\
Our scheme, as previously stated, is to sequentially transform $\mathtt{P}$, using some combination of the objective based approach suggested by (Oa) and (Ob) and the feasible set contraction approach suggested by (Oc), so that it becomes strongly dual. 
Once this transformation is complete, a heuristic for $\mathtt{P}$ is then given by initializing from the corresponding dual solution field and optimizing via some secondary (local) method.   
\subsection{\S Verlan Scheme}
\noindent
Take $\mathtt{P}$ to be a SCQP, with $\bmm{U}$ split into a variable and constant terms $\bmm{U} = \bmm{X}^{-1}_{\bullet}-\bmm{G}_{\circ}$.  
Define $\mathtt{P}\left(\bm{r},\bmm{V}\right)$ to be the program resulting from $\mathtt{P}$ by replacing the linear part of the object, $\bm{s}_{\circ}$ in $f_{o}$ with $\bm{r}$, and the variable operator of the scattering constraint set $\bmm{X}^{-1}_{\bullet}$ with $\bmm{V}^{-1}$. 
Take $\sigma$ and $\gamma$ to be a scrapping parameters and $\epsilon$ to be an expansion parameter. 
In reference to updating program inputs, let $\left(n\right)$, $\left(o\right)$, and $\left(\delta\right)$ superscripts refer to new, old, and difference values respectively, $\bm{x}_{\circledast}^{\left(n\right)} = \bm{x}_{\circledast}^{\left(o\right)} + \bm{x}_{\circledast}^{\left(\delta\right)}$. 
To initialize, set $\bm{r} = \bm{s}_{o}$ and $\bmm{V} = \bmm{X}_{\bullet}$. 
\\ \\
(Pa) \emph{Contract}---Decrease $\bmm{V}\rightarrow 0$ until a few applications of \emph{scrape} (Pc) results in strong duality for $\mathtt{P}\left(\bm{r},\bmm{V}\right)$---practically the condition of strong duality is verified by $\bm{x}_{\circledast}$ being within some tolerance of $\partial F_{\bm{\kappa}}$.  
\\ \\
(Pb) \emph{Expand}---Based on the value of $\epsilon$, let $\bmm{V}$ tend towards to $\bmm{X}_{\bullet}$. 
Solve $\mathtt{D}\left(\mathtt{P}\left(\bm{r},\bmm{V}\right)\right)$. 
If $\bmm{V}$ is sufficient near $\bmm{X}_{\bullet}$, and $\bm{x}_{\circledast}$ is within tolerance of $\partial F_{\bm{\kappa}}$, $\bm{x}_{\circledast}$ can be used to infer a heuristic for $\mathtt{P}$ via a result like eq.~\eqref{dirInfer}.
\\ \\
(Pc) \emph{Scrape}---If $\bm{x}_{\circledast}$ is within tolerance of $\partial F_{\bm{\kappa}}$,  set
$$
  \bm{r}^{\left(n\right)} = \bm{r}^{\left(o\right)} + \sigma\bm{x}_{\circledast}^{\left(o\right)},
$$ 
and return to \emph{expand} (Pb). 
If $\bm{x}_{\circledast}$ is not within tolerance of $\partial F_{\bm{\kappa}}$ set 
$
  \bm{r}^{\left(n\right)} = \bm{r}^{\left(o\right)} + \gamma\bm{x}_{\circledast}^{\left(o\right)}.
$
Continue applying \emph{scrape} until either (a): $\bm{x}_{\circledast}^{\left(n\right)}$ is sufficiently near $\partial F_{\bm{\kappa}}$, or (b): a limit number of \emph{scrape}s is reached.  
If exiting under (a), proceed to a new application of \emph{expand} (Pb). 
If exiting under (b), reduce $\epsilon$, revert $\mathtt{P}\left(\bm{r},\bmm{V}\right)$ to the last state with strong duality, and then apply \emph{expand}~\footnote{It may also be sensible to normalize $\bm{r}^{\left(n\right)}$ or include some removal of $\bm{r}^{\left(o\right)}$ in the update, i.e. $\bm{r}^{\left(n\right)} = \left(1-\gamma\right)\bm{r}^{\left(o\right)} + \gamma\bm{x}_{\circledast}^{\left(o\right)}$.}.
\\ 
In carrying out the \emph{expand} and \emph{scrape} steps, with idealized precision and control, the algorithm tracks the evolution of a feasible solution near the maximum of the original objective function. 
The $\sigma$ parameter is used to mediate between this goal and the computational expedient of maintaining a tolerable rate of expansion towards the original constraint set. 
$\sigma > 0$ causes the $\bm{r}$ to continuously shift towards the current solution, improving the chances that strong duality will be preserved in the next application of \emph{expand}, which will presumably improve the numeric conditioning of the dual program.  
By setting $\gamma = 0$, deviation between $\bm{r}$ and $\bm{s}_{o}$ is forestalled for as long as possible. 
We conjecture that smaller values of $\sigma$ and $\gamma$ will generally lead to greater fidelity between $\bm{r}$ and $\bm{s}_{o}$ at the cost of more stringent requirements on numerical precision---some \emph{expand} step will eventually require prolonged scraping, and scraping will be slower. 
\\ \\
As suggested by (Ob), \emph{scrape} (Pc) is a method of pushing $\bm{x}_{\circledast}$ towards $\partial F_{\bm{\kappa}}$. 
Suppose $\bm{x}\in S_{\mathtt{P}}\Rightarrow\lVert\bm{s}_{o}\rVert\geq\lVert\bmm{A}_{o}\bm{x}\rVert$: the maximum of $f_{o}$ is in $\partial S_{\mathtt{P}}$~\footnote{If it is not clear that this assumption holds, an additional normalization step can be added to \emph{scrape}.}.
Because $f_{o}\left(\bm{r}^{\left(n\right)},\bm{x}_{\circledast}^{\left(n\right)}\right)\geq f_{o}\left(\bm{r}^{\left(n\right)},\bm{x}_{\circledast}^{\left(o\right)}\right)$ and $f_{o}\left(\bm{r}^{\left(o\right)},\bm{x}_{\circledast}^{\left(o\right)}\right)\geq f_{o}\left(\bm{r}^{\left(o\right)},\bm{x}_{\circledast}^{\left(n\right)}\right)$, $\Re\left[\bm{y}^{\left(\delta\right)\dagger}\bm{x}_{\circledast}^{\left(\delta\right)}\right]\geq 0$.
Hence,
\begin{align}
  \lVert\bm{x}_{\circledast}^{\left(n\right)}\rVert^{2} &= \lVert\bm{x}_{\circledast}^{\left(o\right)}\rVert^{2} + 2\Re\left[\bm{x}_{\circledast}^{\left(o\right)\dagger}\bm{x}_{\circledast}^{\left(\delta\right)}\right] + \lVert\bm{x}_{\circledast}^{\left(\delta\right)}\rVert^{2} 
  \nonumber \\
  &= \lVert\bm{x}_{\circledast}^{\left(o\right)}\rVert^{2} + \frac{2}{\gamma}\Re\left[\bm{y}^{\left(\delta\right)\dagger}\bm{x}_{\circledast}^{\left(\delta\right)}\right] + \lVert\bm{x}_{\circledast}^{\left(\delta\right)}\rVert^{2} \geq \lVert\bm{x}_{\circledast}^{\left(o\right)}\rVert^{2}.
  \nonumber
\end{align}
Directly, the scrape protocol guarantees that the magnitude of $\bm{x}_{\circledast}$ never decreases. 
So long as $\bm{x}_{\circledast}$ is not trapped at a high curvature point of the Sion set in the interior of $F_{\bm{\kappa}}$, it will eventually push $\bm{x}_{\circledast}$ to $\partial F_{\bm{\kappa}}$. 
This potential failure mode of \emph{scrape}, if employed without the other two protocols, can be understood from imagining a small alteration to fig.~2(a). 
If the maximal bulge in the Sion set (purple region) were to occur below the initial objective direction (first red arrow), rather than to the left, then scrapping would cause $\bm{s}_{o}$ to converge to the bulge direction. 
Under contraction, $S_{\mathtt{P}}$ increasingly coincides with $C_{\mathtt{P}}$, and on $C_{\mathtt{P}}$ such a possibility is effectively ruled out.
\\ \\
We have not tested enough scenarios to provide useful rules of thumb for how to best combine protocols (Pa)--(Pc) in a typical application, c.f. ref.~\cite{chao2025blueprints}, and are unaware of any rigorous optimality. 
Nevertheless, the definition of $\mathtt{S}\left(\mathtt{P}\left(\bm{r},\bmm{V}\right)\right)$ implies that $\bm{x}_{\circledast}$ is usually a point in $S_{\mathtt{P}}\left(\bmm{V}\right)$ with a near optimal real projection along $\bm{r}$---if $f_{o}$ is linear, then $\bm{x}_{\circledast}$ is the point in $S_{\mathtt{P}}\left(\bmm{V}\right)$ with the largest $\bm{r}$ coefficient. 
When $\bm{r}$ must be altered to induce strong duality, absent additional information, shifting towards $\bm{x}_{\circledast}$ is at least rationally motivated. 
Much like following the gradient, the basic idea is that combining a series of (approximately) optimal local choices generally creates a good heuristic.
\subsection{\S Outlook}
\noindent 
While the success of the example protocol presented in ref.~\cite{chao2025blueprints} highlights the immediate applicability of verlan approaches to contemporary photonic device design problems, there are several outstanding computational challenges that must be met before the full potential of these methods can be realized.
At the same time, there are also a number of promising directions for further advancement that can be explored immediately. 
\\ \\
Through the imposition of spherical symmetries, and restrictions on the number of addressable constraints, duality relaxations have already been used to compute performance limits for design volumes exceeding one thousand cubic wavelengths in size~\cite{molesky2020t,molesky2020hierarchical}.
Nevertheless, limit programs for three-dimensional photonic inverse design problems are presently limited by computational costs. 
In either the SDP or dual framework, solving a limit program requires a large number of verifications that an overall system matrix (of a size proportional to the discretization) is positive semi-definite. 
Although sparse (chordal) Cholesky factorization~\cite{chen2008algorithm,rennich2014accelerating,cao2022framework} offers important speed-ups compared to calculation of the smallest eigenvalue for performing this check, memory requirements for these techniques are imposing for wavelength-scale three-dimensional systems, and checking positive-definiteness unequivocally remains the critical bottleneck. 
Further, even if this bottleneck is addressed, there is still a need to perform hundreds to thousands of inverse linear solves on poorly conditioned systems, which become increasingly difficult as the system size grows.
These issues do not need to be completely solved before verlan approaches could have substantial utility in large-scale photonic inverse design---e.g. based the performance of the global curves shown in ref.~\cite{chao2025blueprints}, the idea of using of weakly constrained QCQPs to carry out ``relaxed'' verlan inverse design over patches spanning tens of wavelengths in order to improve the performance of large scale meta-lenses is well motivated and within the reach of current codebases---but to perform ``full'' verlan inverse design over three dimensional volumes surpassing  one hundred cubic wavelengths, without the use of simplifying symmetries, significant developments will need to be made in the numerics of both limit programs and photonic solvers~\cite{xue2023fullwave,mao2024towards,pestourie2025fast}. 
\\ \\
Bearing these limitations in mind, there are also many avenues for progress that can be explored with existing tools. 
For example, barring a handful of exceptions~\cite{zhang2021conservationlawbased,wisal2024optimal,mishra2024classically}, the QCQP--convex relaxation procedure has been applied to photonics almost exclusively in the frequency domain. 
Eliding the the nuances of properly describing material properties in the time domain~\cite{hordt1998calculation}, this recipe, and by extension verlan inverse design, applies essentially without modification in form to the time domain~\cite{liska2024excitation}. 
As the few references given above indicate, it is likely that, beyond determining device structures for scenarios in which QCQP performance bounds have been previously shown, verlan approaches can be productively applied to dynamic processes, enabling a variety of applications that have yet to be considered~\cite{hammond2022high,gedeon2023time,garg2024inverse}.
\\ \\
Closer to the content of this manuscript, there is also probably much to be gained from integrating our verlan scheme with the approximate rank penalization methods typically used in SDP relaxations~\cite{fazel2002matrix,luo2010semidefinite,madani2014convex,gertler2025many}. 
Results concerning the relation between the rank of SDP solutions and maximal clique size in the imposed constraints seem likely to be connected (in some way) to deviations of the Sion set from the convex hull of feasible points~\cite{vandenberghe2015chordal}. 
Additionally, both approximate rank penalization and the \emph{scrape}-\emph{expand} protocols we have proposed have the primary goal of transforming the underlying the QCQP, ideally through an optimally small modification, so that it becomes strongly dual. 
Given that SDP and duality convexifications are related by duality, and in fact meet the requirements of Sion's theorem in the case of a scattering constraint set, it is clear that either procedure can be translated from one form of convex relaxation to the other. 
Examining the resulting interaction appears to us as a promising prospect for improving understanding, and ultimately performance, of general inverse design heuristics. 
\section{Appendices}
\noindent 
The following four sections present some previously established results supplementing the main text.
\subsection{\S Continuity Lemma}
\noindent
\emph{Lemma A.} If $f:X\times C\rightarrow\mathbb{R}$ is continuous, and $C$ is compact, then $g:X\rightarrow\mathbb{R},~x\mapsto\sup_{c\in C} f\left(x,c\right)$ and $h:X\rightarrow\mathbb{R},~x\mapsto\inf_{c\in C} f\left(x,c\right)$ are continuous. 
\begin{proof}
  Because continuity implies point continuity, about each $\left(x,c\right)\in X$, $\forall \delta > 0$, there is a product neighbourhood $N_{c}^{\left(x,c\right)}\times M_{c}^{\left(x,c\right)}$ such that $f\left[N_{x}^{\left(x,c\right)}\times M_{c}^{\left(x,c\right)}\right]\subseteq\left(f\left(x,c\right) - \delta, f\left(x,c\right) + \delta\right)$. 
  A finite collection of $M_{c}^{\left(x,c\right)}$ subsets, indexed by $F_{x}=\left\{c_{i}\right\}$, cover $C$ for each $x\in X$.
  Take $N_{x} = \cap_{c_{i}\in F_{x}}N_{x}^{\left(x,c_{i}\right)}$. 
  $\forall c\in C~f\left[N_{x}\times\left\{c\right\}\right]\subseteq\left(f\left(x,c\right) - \delta, f\left(x,c\right) + \delta\right)$.
  Via the product topology, $f_{x}:C\rightarrow\mathbb{R},~c\mapsto f\left(x,c\right)$, is continuous. 
  Hence, for any given $x\in X$, there is a $c_{u}\in C \ni g\left(x\right) = f\left(x,c_{u}\right)$. 
  As such, $g\left[N_{x}\right]\subseteq\left(g\left(x\right) - \delta, g\left(x\right) + \delta\right)$, and $g$ is continuous. 
  For any continuous $f:X\times C\rightarrow\mathbb{R}$, $\inf_{c\in C} f\left(x,c\right) = \sup_{c\in C} -f\left(x,c\right)$. 
\end{proof}
\subsection{\S Schur Complement}
\noindent
Recall that if a block matrix is either upper or lower diagonal then its determinant is given by the product of the determinants of its main diagonal blocks. 
So long as $\bmm{A}$ is invertible, the expansion 
\begin{align}
  &\begin{bmatrix}
    \bmm{A} & \bmm{C} \\
    \bmm{B} & \bmm{D}
  \end{bmatrix} = 
  \begin{bmatrix}
    \bmm{I} & \bmm{0} \\
    \bmm{B}\bmm{A}^{-1} & \bmm{I}
  \end{bmatrix}
  \begin{bmatrix}
    \bmm{A} & \bmm{C} \\
    \bmm{0} & \bmm{D}-\bmm{B}\bmm{A}^{-1}\bmm{C}
  \end{bmatrix} =
  \label{Schur} \\
  &\begin{bmatrix}
    \bmm{I} & \bmm{0} \\
    \bmm{B}\bmm{A}^{-1} & \bmm{I}
  \end{bmatrix}
  \begin{bmatrix}
    \bmm{A} & \bmm{0} \\
    \bmm{0} & \bmm{D}-\bmm{B}\bmm{A}^{-1}\bmm{C}
  \end{bmatrix}
  \begin{bmatrix}
    \bmm{I} & \bmm{A}^{-1}\bmm{C} \\
    \bmm{0} & \bmm{I}
  \end{bmatrix}
  \nonumber
\end{align}
indicates that the total matrix is positive semi-definite iff $\bmm{A}\succ 0$ and $\bmm{D}-\bmm{B}\bmm{A}^{-1}\bmm{C}\succ 0$. 
Due to this critical property, $\bmm{D} - \bmm{B}\bmm{A}^{-1}\bmm{C}$ is referred to as the Schur complement of $\bmm{A}$. 
An analogous expansion can be carried out if $\bmm{D}$ is invertible, showing that the block matrix is positive definite iff $\bmm{D}\succ 0$ and $\bmm{A}-\bmm{C}\bmm{D}^{-1}\bmm{B}\succ 0$.
\subsection{\S Equivalence with Semi-Definite Relaxation}
\noindent
The dual, in the wider setting of QCQPs, is equivalent to the semi-definite program relaxation. 
Taking $f_{o}\left(\bm{x}\right)$ in eq.~\eqref{SCQP} to be a quadratic function, this equivalence is realized as follows. 
Using the correspondence between bilinear forms and tensor products, any quadratic constraint $f_{j}\left(\bm{x}\right)\geq 0$ is equivalent to 
\begin{equation}
  \tr{
  \begin{bmatrix}
    -\bmm{A}_{j} & \bm{s}_{j} \\ 
    \bm{s}_{j}^{\dagger} & v_{j}
  \end{bmatrix}
  \begin{bmatrix}
    \bm{x}\bm{x}^{\dagger} & \bm{x} \\ 
    \bm{x}^{\dagger} & 1
  \end{bmatrix}}
  \geq 0.
\end{equation}
Letting $\bmm{H}_{k}$ denote the block matrix associated with a quadratic function in the above manner, by increasing the dimension of $\bm{x}$ by one, $\bm{x}\rightarrow\tilde{\bm{x}}$, any QCQP $\mathtt{P}$ can be placed in the homogeneous form
\begin{align}
  &\max_{\tilde{\bm{x}}\in\mathbb{C}^{n+1}} \tr{\bmm{H}_{o}\bmm{Y}}
  \label{homoForm} \\
  &\ni\left(\forall j\in J\right)~\tr{\bmm{H}_{j}\bmm{Y}}\geq 0
  \nonumber \\
  &\land \tilde{x}_{n+1}^{2} = 1,
  \nonumber
\end{align}
where $\bmm{Y} = \tilde{\bm{x}}\tilde{\bm{x}}^{\dagger}$---a dummy variable has been used to increase the dimension of the problem by one, if $\tilde{\bm{x}}$ is a solution of eq.~\eqref{homoForm}, then $\bm{x} = \left(\tilde{x}_{1}/\tilde{x}_{n+1},\tilde{x}_{2}/\tilde{x}_{n+1},\ldots,1\right)$ is a solution of the original QCQP.
The dual of eq.~\eqref{homoForm}, using $\bm{x}_{\circledast} = \bmm{A}_{\bm{\psi}}^{-1}\bm{s}_{\bm{\psi}}$ and $\bm{\psi}$ to denote the inclusion of the objective in the $\bm{\phi}$ sum notation used in the main text, is 
\begin{equation}
  \inf_{\bm{\phi}\in\Psi_{\mathtt{P}}}\bm{s}_{\bm{\psi}}^{\dagger}\bmm{A}^{-1}_{\bm{\psi}}\bm{s}_{\bm{\psi}} + v_{\bm{\psi}}.
  \nonumber
\end{equation}
Via the Schur complement construction, the dual of a QCQP is thus equivalent to 
\begin{align}
  \inf_{\alpha\in\mathbb{R},\bm{\phi}\in\mathbb{H}^{J}}\alpha
  \ni
  \begin{bmatrix}
    -\bmm{A}_{\bm{\psi}} & \bm{s}_{\bm{\psi}} \\
    \bm{s}_{\bm{\psi}}^{\dagger} & v_{\bm{\psi}} - \alpha 
  \end{bmatrix}\preceq 0.
  \label{SchurDual}
\end{align}
The constraint that the block matrix appearing in eq.~\eqref{SchurDual}, $\bmm{H}^{\alpha}_{\bm{\psi}}$, is negative semi-definite is equivalent to checking that $\bmm{b}_{l}^{\dagger}\bmm{H}^{\alpha}_{\bm{\psi}}\bmm{b}_{l}\leq 0$ for a collection of basis vectors $\left\{\bmm{b}_{l}\right\}_{l\in L}$ spanning $\mathbb{C}^{n+1}$. 
Since the sum of a collection of such constraints can be represented by a positive semi-definite matrix $\bmm{B}$, the Lagrangian of eq.~\eqref{SchurDual} is 
$$
  \mathcal{L}\left(\alpha,\bm{\phi},\bmm{B}\right) = 
  \alpha \left(1 - \bmm{B}_{n+1,n+1}\right) + \tr{\bmm{H}_{\bm{\psi}}\bmm{B}}. 
$$
If either $\bmm{B}_{n+1,n+1} \neq 1$ or $\tr{\bmm{H}_{j}\bmm{B}} < 0$ for some $j\in J$, then the infimum of this form over $\alpha$ and $\bm{\phi}\in\mathbb{H}^{J}$ is $-\infty$; if $\bmm{B}_{n+1,n+1} = 1$ and $\left(\forall j\in J\right)~\tr{\bmm{H}_{j}\bmm{B}}\geq 0$, then the infimum is $\tr{\bmm{H}_{o}\bmm{B}}$. 
The dual of eq.~\eqref{SchurDual}, accordingly, is 
\begin{align}
  &\max_{\bmm{B}\in\mathbb{C}^{\left(n+1\right)^{2}},\bmm{B}\succeq 0}~\tr{\bmm{H}_{o}\bmm{B}}
  \label{dualDual} \\
  &\ni\left(\forall j\in J\right)~\tr{\bmm{H}_{j}\bmm{B}}\geq 0
  \land \bmm{B}_{n+1,n+1} = 1.
  \nonumber
\end{align}
eq.~\eqref{dualDual} is equivalent to eq.~\eqref{homoForm} under the semi-definite relaxation that $\bmm{Y}$ is simply some positive semi-definite matrix---simply forgetting that the rank of $\bmm{Y}$ as written in eq.~\eqref{homoForm} is one. 
As can be shown in greater generality~\cite{arnol2013mathematical}, the duality relaxation is an involution, and taking the dual of \eqref{dualDual}.
The Lagrangian of eq.~\eqref{dualDual} is
$$
  \mathcal{L}\left(\alpha,\bm{\phi},\bmm{Z},\bmm{B}\right) = 
  \tr{\left(\bmm{H}_{\bm{\psi}}+\bmm{Z}\right)\bmm{B}} + \alpha\left(1-\bmm{B}_{n+1,n+1}\right). 
$$
Since $\bmm{Z}\succeq 0$, its addition can never lead to a smaller maximum over $\bmm{B}$, and so minimizing the dual always implies $\bmm{Z} = \bmm{0}$. 
Avoiding infinite values requires that $\bmm{H}_{\bm{\psi}}^{\alpha}\preceq 0$, and under this condition $\bmm{B}_{\circledast} = \bmm{0}$. 
Making these substitutions, the dual returns to eq.~\eqref{SchurDual}. 
\\ \\
One way of interpreting the SDP relaxation, from the perspective of Lagrange duality, is that conversion of the dual into a SDP effectively introduces additional copies of the dual in proportion to the rank of the solution matrix: each vector in the eigen expansion of $\bmm{B}$ defines a multiplier function that is equivalent to the dual up to a lack of implicit dependence between the vector and the multiplier values. 
By summing these ``dual'' copies, $\bmm{H}_{\bm{\psi}}^{\alpha}$ is confined to the semi-definite cone, so that the aggregate system becomes strongly dual---some of the individual fields may violate constraints, but as a sum every imposed constraint is satisfied. 
If only a single copy of the dual is required, then there is an immediate correspondence between the two descriptions and the computed solutions are identical. 
If the rank of the SDP is greater than one (i.e. multiple copies are required), an additional (convex) optimization is needed to relate the two solutions. 
\subsection{\S QCQP Transformation for Non-Local Scattering}
\noindent
Regardless of whether $\bmm{X}$ is spatially local (diagonal in the spatial basis) or not, the scattering relation developed in \emph{\S QCQPs in Photonic Device Design} is equivalently given in terms of polarization distributions as 
\begin{equation}
  \left(\bmm{I}-\bmm{X}\bmm{G}_{\circ}\right)\bm{j}_{s} = \bmm{X}\bmm{G}_{\circ}\bm{j}_{i}.
  \label{sctCon}
\end{equation}
Supposing that possible instances of $\bmm{X}$ are controlled by a set of design parameters $\theta_{j}\in\left\{0,1\right\}$~\footnote{As in the main text, the quadratic equations resulting from eq.~\eqref{sctCon} can be adapted to continuous design parameters by replacing equalities with inequalities.}, let $\bmm{Q}_{j}$ be a witness of the $j^{th}$ design parameter in the sense that $\bmm{Q}_{j}\bmm{X} = \bmm{0}$ whenever $\theta_{j} = 0$ and $\bmm{Q}_{j}\bmm{X} = \bmm{R}_{j}$, for some fixed operator $\bmm{R}_{j}$, whenever $\theta_{j} = 1$.
Applying $\bmm{Q}_{j}$ to \eqref{sctCon}, it follows that $\bmm{Q}_{j}\bm{j}_{s} = \bm{0}$ whenever $\theta_{j} = 0$. 
Letting $\bmm{X}_{\bullet}$ signify the scattering potential given by $\left(\forall j\in J\right)~\theta_{j} = 1$,
\begin{equation}
  \bm{j}_{s}^{\dagger}\bmm{Q}_{j}^{\dagger}\bmm{P}\bmm{Q}_{j}\left(\bmm{I}-\bmm{X}_{\bullet}\bmm{G}_{\circ}\right)\bm{j}_{s} = \bm{j}_{s}^{\dagger}\bmm{Q}_{j}^{\dagger}\bmm{P}\bmm{Q}_{j}\bmm{X}_{\bullet}\bmm{G}_{\circ}\bm{j}_{i}. 
  \label{allOn}
\end{equation}
is true for any linear operator $\bmm{P}$.
As in eq.~\eqref{sctCnt}, the influence of the design in eq.~\eqref{allOn} has been shifted from the scattering potential $\bmm{X}$ to the scattered current $\bm{j}_{s}$, allowing eq.~\eqref{sctCon} to be restated as a collection of scattering constraints. 
Such a procedure can be implemented whenever $\bmm{X} = \sum_{j}\theta_{j}\bmm{D}_{j}$ is ``image separable'': to each $j$ there is a linear functional $\bm{l}_{j}^{\dagger}$ such that $\bm{l}_{j}^{\dagger}\bmm{D}_{j}\neq\bmm{0}^{\dagger}$ and $\bm{l}_{j}^{\dagger}\bmm{D}_{k}=\bmm{0}^{\dagger}$ for all $k\neq j$. 
\bibliography{esaLibF}
\end{document}